\documentclass[10pt,leqno]{amsart}
\usepackage{verbatim}
\usepackage{enumerate}
\usepackage[mathscr]{eucal}
                                           
\title{The coarea formula for Sobolev mappings}
\author{Jan Mal\'y, David Swanson and William P.\ Ziemer}
\address{Faculty of Mathematics and Physics, Charles University --
KMA, Sokolovsk\'a 83, 18675 Praha 8, Czech Republic}
\email{maly@karlin.mff.cuni.cz}
\address{Department of Mathematics, Texas A{\&}M University,
College Station, TX, 77843}
\email{dswanson@math.tamu.edu}
\address{Department of Mathematics, Indiana University,
Bloomington, IN, 47405}
\email{ziemer@indiana.edu}
\thanks{The research of the first author is supported in part by the  
Research Project MSM 113200007 from the Czech Ministry of Education,
Grant No. 201/00/0767 from the Grant Agency of the Czech republic
(GA\,\v CR) and Grant No. 165/99 from the Grant Agency 
of Charles University (GA\,UK)}

\numberwithin{equation}{section}

\swapnumbers
\newtheorem{theorem}{Theorem}[section]
\newtheorem{lemma}[theorem]{Lemma}
\newtheorem{prop}[theorem]{Proposition}
\newtheorem{cor}[theorem]{Corollary}
\theoremstyle{definition}

\newtheorem{rem}[theorem]{Remark}

\newtoks\by
\newtoks\paper
\newtoks\book
\newtoks\jour
\newtoks\yr
\newtoks\pages
\newtoks\vol
\newtoks\publ
\def\ota{{\hbox\vol{???}}}
\def\cLear{\by=\ota\paper=\ota\book=\ota\jour=\ota\yr=\ota
\pages=\ota\vol=\ota\publ=\ota}
\def\endpaper{\the\by, {\the\paper},
\textit{\the\jour} \textbf{\the\vol} (\the\yr), \the\pages.\cLear}
\def\endbook{\the\by, \textit{\the\book}, \the\publ.\cLear}
\def\endprep{\the\by, \textit{\the\paper}, \the\jour.\cLear}
\def\endproc{\the\by, \textit{\the\paper}, \the\book, \the\publ,
\the\yr, \the\pages.\cLear}
\def\name#1#2{#1, #2}

\newcommand{\f}{\bar f}
\newcommand{\loc}{\mathrm{loc}}
\newcommand{\diam}{\mathrm{diam}}

\newcommand{\eps}{\varepsilon}
\newcommand{\R}{\mathbb{R}}

\newcommand{\hm}{\mathcal{H}}

\newdimen\vintkern
\def\cp{\gamma }
\def\en{\mathbb N}
\def\er{\mathbb R}
\def\ef{\mathbf F}
\def\erm{\er^m}
\def\G{\mathcal G}
\def\lip{\operatorname{Lip}}
\def\unitball{\boldsymbol\alpha}

\def\L{\mathcal L}
\def\grass{\mathbf{G}}
\def\vint{-\kern\vintkern\int}
\def\eqn#1$$#2$${\begin{equation}\label#1#2\end{equation}}
\def\intav#1{\vint_{#1}}

\begin{document}
\begin{abstract}
We extend Federer's coarea formula to mappings $f$ belonging 
to the Sobolev class $W^{1,p}(\R^n;\R^m)$, $1 \le m < n$, $p>m$, and more
generally, to mappings with gradient in the Lorentz space $L^{m,1}(\R^n)$.
This is accomplished by showing that the graph of $f$ in $\R^{n+m}$
is a Hausdorff $n$-rectifiable set.
\end{abstract}

\maketitle
\section{Introduction}

The following equality, first proved by Federer in \cite{F3} and termed
the coarea formula, has proved to be a fundamental tool in analysis 
with a wide variety of applications:
\begin{equation}
\label{coarea1}
\int_{\Omega} g(x) |J_mf(x)| \, dx = \int_{\R^m} \int_{f^{-1}(y)}
g(x) \, d\hm^{n-m}(x) \, dy,
\end{equation}
where $\Omega\subset\R^n$ is an open set, $f : \Omega \to \R^m$ is Lipschitz, 
$J_mf$ is its $m$-dimensional Jacobian,
$g:\Omega \to \R$ is integrable, and $1 \le m < n$.  Recall that
$|J_{m}f|$ is the square root of the sum of the squares of
the determinants of the $m$ by $m$ minors of the differential of $f$.
Using standard approximation methods (\ref{coarea1}) may be obtained
from the special case where $g$ is the indicator of a measurable
set $E \subset \R^n$:
\begin{equation}
\label{coarea2}
\int_E |J_mf(x)| \, dx = \int_{\R^m} \hm^{n-m}(E \cap f^{-1}(y)) \, dy.
\end{equation}

Because of the usefulness of the coarea formula, a natural and
compelling question is whether it will admit an extension to a larger
class of mappings, in particular, to the class 
$W^{1,1}_{\loc}(\Omega;\R^m)$ of Sobolev mappings.
When $m=1$, it was shown by Federer \cite{GMT} that \eqref{coarea2}
is valid for mappings $f \in W^{1,p}(\Omega)$, $p \ge 1$,
provided that $f$ is precisely represented, see 
\eqref{defpreciserep} below.  The case $p > n$, including $m \ge 2$, 
has been considered by Van der Putten \cite{VP}. 

One of the main objectives of this paper is to establish (\ref{coarea2})
for mappings $f$ belonging to the Sobolev class $W^{1,p}_{\loc}(\Omega;\R^m)$.
The following result is new for $m > 1$ and the proof simplifies the
existing proof for $m=1$.

\begin{theorem}
\label{t:main}
Suppose that $1\le m \le n$ and that 
$f \in W_{\loc}^{1,p}(\Omega;\R^m)$ is precisely represented, 
where either $p > m$ or $p \ge m = 1$.
Then $f^{-1}(y)$ is countably $\mathcal{H}^{n-m}$ rectifiable
for almost all $y \in \R^m$ and the coarea formula
\eqref{coarea2} holds
for all measurable sets $E \subset \Omega$.
\end{theorem}

In the following theorem
we denote by $\bar f : \Omega \to \R^{n+m}$
the graph mapping $\bar f (x) = (x,f(x))$ and by $\mathcal{G}_f$ the
graph of $f$:  $\mathcal{G}_f = \bar f(\Omega) \subset \R^{n+m}$.

\begin{theorem}
\label{t:graph}
Suppose that $1\le m \le n$ and that 
$f \in W_{\loc}^{1,p}(\Omega;\R^m)$ is precisely represented, 
where either $p > m$ or $p \ge m = 1$.
Then $\G_f$ is countably $\mathcal{H}^n$ rectifiable and 
\eqn{area2}
$$
\hm^n(\bar f(E)) = \int_E |J_n \bar f(x)| \, dx
$$
for all measurable sets $E \subset \Omega$. In particular, $\hm^n(\bar
f(E))=0$ whenever $\L^{n}(E)=0$.
\end{theorem}

The hypothesis $p > m$ is necessary. In \cite{C}, \cite{MM}, an example
of a continuous mapping $f \in W^{1,m}(\R^m;\R^m)$ is exhibited which 
has an almost everywhere vanishing Jacobian and maps
an interval $I$ onto an $m$-cube. This may be easily modified in the case
$m \ge 2$ to an example of a continuous mapping 
$f \in W_{\loc}^{1,m}(\R^n;\R^m)$
which has an almost everywhere vanishing $m$-Jacobian 
and which maps each set of the form $I\times \R^{n-m}$
onto an $m$-cube (cf.\ \cite{Ha}), therefore violating the conclusions
of Theorems~\ref{t:main} and~\ref{t:graph}. 
 However, if we turn
our attention to mappings whose
gradients belong to the finer scale of Lorentz
spaces, we obtain the following extension of the above
stated theorems.

\begin{theorem}
\label{t:lorentz}
Suppose that $1\le m\le n$, that 
$f \in W_{\loc}^{1,1}(\Omega;\R^m)$ 
is precisely represented, and that $|\nabla f|\in L^{m,1}(\Omega)$.
Then $f^{-1}(y)$ is countably $\mathcal{H}^{n-m}$ rectifiable
for almost all $y \in \R^m$,
the graph $\G_f$ is countably $\mathcal{H}^n$ rectifiable and 
the coarea formula \eqref{coarea2} and area formula \eqref{area2}
hold
for all measurable sets $E \subset \Omega$. 
\end{theorem}

The case $m=n$ of Theorem \ref{t:lorentz} (and even the area 
counterpart $m \ge n$) follows from the work
of Kauhanen, Koskela, and Mal\'y \cite{KKM} who generalized the
area formula of Marcus and Mizel \cite{MMi} to functions whose
gradients lie in $L^{n,1}(\Omega)$.

It is also possible to obtain results dealing with the borderline
case $p=m$ for mappings which are H\"{o}lder continuous.  The
following theorem provides a coarea
counterpart to the results in \cite{MM} on the area
formula for H\"older continuous mappings in $W^{1,n}$
.
(See also the recent developments in 
\cite{M-area}, \cite{FMy}, \cite{M-ag}).

\begin{theorem}\label{t:hoelder}
Suppose that $1 \le m \le n$ and that
$f\in W^{1,m}(\Omega;\erm)$ is a H\"older continuous mapping.
Then $f^{-1}(y)$ is countably $\mathcal{H}^{n-m}$ rectifiable
for almost all $y \in \R^m$ and  
\[
\int_E |J_m f(x)| \, dx = \int_{\R^m} \hm^{n-m}(E \cap f^{-1}(y)) \, dy
\]
for all measurable sets $E \subset \Omega$.
\end{theorem}

Theorem~\ref{t:hoelder} is proven independently of any area estimates on the
$\mathcal{G}_f$ and in fact the proof can be modified to give an independent
proof of Theorem~\ref{t:main} in the case $p > m$, cf.\
Remark~\ref{r:proof} below. 

Weaker variants of the coarea and area formulas use the integralgeometric
measure
 $\mathcal{I}^q$ instead of Hausdorff measure $\mathcal{H}^q$. 
Haj{\l}asz \cite{Ha} proved such a version of 
the coarea formula for $W^{1,p}$-mappings ($p>m$).
We present a version of Theorem~\ref{t:graph} 
in the borderline case $p=m$ under the assumption of
H\"{o}lder continuity. We have not been able to obtain such a result
with the Hausdorff measure.
  
\begin{theorem}
\label{t:hoelderarea}
Suppose that $1 \le m \le n$ and that $f \in W^{1,m}(\Omega;\R^m)$
is a H\"{o}lder continuous mapping.  Then $\mathcal{G}_f$ is 
countably $\mathcal{I}^n$ rectifiable and
\[
\mathcal{I}^n (\bar f(E)) = \int_E |J_n \bar f(x)| \, dx
\]
for all measurable sets $E \subset \Omega$.
\end{theorem}

The results listed above are easily localized from $\R^n$ to $\Omega$. 
For the proofs we may assume that $\Omega=\R^n$ and that the global 
norms of the considered mappings are finite.
  Since our results follow
from the area formulas in \cite{MMi} and \cite{KKM} when $m=n$ we also
assume throughout that $m<n$.

\section{Preliminaries}

We denote by $\L^{n}(E)$ the Lebesgue measure of a set $E \subset \R^n$
and by $\hm^q(E)$ the $q$-dimensional Hausdorff measure of $E$.  
If $q$ is an integer then $E$ is said to be countably 
$\hm^q$ rectifiable
if there exist subsets $E_k \subset \R^q$ and Lipschitz mappings
$g_k : E_k \to \R^n$ with the property that
\[
\hm^q \left( E \setminus \bigcup_{k=1}^\infty  
g_k(E_k) \right) = 0.
\]

The $q$-dimensional Hausdorff content $\hm^q_\infty(E)$ of $E$ 
is the infimum of the sums 
\[ 
\sum_{j=1}^\infty \unitball_q\left( \frac{\diam\, E_j}{2} \right)^q 
\] 
corresponding to all countable coverings of $E$ by sets 
$\{E_j\}_{j=1}^\infty$. 
Here 
$$ 
\unitball_q = \frac{\pi^{q/2}}{\Gamma(\frac q2+1)}. 
$$ 
which is the volume of the unit ball in $\er^q$ if $q$ is integer. 
Observe that $\hm^q(E) = 0$ if and 
only if $\hm^q_\infty(E) = 0$. 
The integral average
of a function $g$ over the set $E$ is defined as
\[
\intav{E} g(y) \, dy = \frac{1}{\L^n(E)} \int_E g(y) \, dy.
\]

Given $1 \le p \le \infty$
the Sobolev space $W^{1,p}(\Omega)$ consists of those functions
$f \in L^p(\Omega)$ with the property that the distributional
gradient $Df$ of $f$ may be identified with a function $\nabla
f \in L^p(\Omega;\R^n)$. 
$W^{1,p}(\Omega)$ is a Banach space with respect to the norm
\[
\|f\|_{1,p,\Omega} := \|f\|_{L^p(\Omega)} + \| \nabla f\|_{L^p(\Omega)}.
\]
The class $W^{1,p}(\Omega;\R^m)$, $m \ge 1$, consists of those
mappings $f : \Omega \to \R^m$ whose component functions each
belong to $W^{1,p}(\Omega)$. Qualitative properties like 
Lebesgue points, approximate differentiability, etc.\ may be
investigated componentwise.

We will also consider a refined scale of Sobolev spaces consisting
of those
functions with gradients in Lorentz and Orlicz spaces.
If $g$ is a measurable function on $\Omega$, we define the \emph{distribution
function} of $g$ as
$$
\mu_g(s) = \L^n\bigl(\{x\in\Omega\colon |g(x)|>s\}\bigr).
$$
We set
$$
\|g\|_{L^{m,1}(\Omega)} = \int_0^{\infty}\mu_g(s)^{1/m}\,ds.
$$
We say that $g$ belongs to the \emph{Lorentz space} $L^{m,1}(\Omega)$ if
$\|g\|_{L^{m,1}(\Omega)}<\infty$. 

It is known that this space can be expressed as a union of Orlicz spaces.
A function $\ef\colon [0,+\infty)\to\er$ 
is said to be a \emph{Young function}
if $\ef$ is convex, nonnegative and satisfies
$$
\ef(t)=0\iff t=0.
$$
The Orlicz space $L^{\ef}(\Omega)$ is then defined as 
the space of all measurable functions $g$ on $\Omega$ for 
which there exists $\lambda>0$ such that
$$
\int_{\Omega}\ef(|g|/\lambda)\,dx\le 1.
$$

The following result is proven in \cite{MSZ}.

\begin{prop}
\label{p:msz}
Let $g\in L^{m,1}(\Omega)$ and $1<q<m<p$. Then
there is a $\mathcal C^1$ Young function $\ef$ and a constant
$C = C(m,p,q) > 0$ such that
\eqn{ef1}
$$
\int_{\Omega}\ef(|g|)\,dx
\le C\|g\|_{L^{m,1}(\Omega)}^m
$$
\eqn{ef2}
$$
\int_{0}^{\infty}\ef'(t)^{-\frac1{m-1}} \,dt \le C
$$
and
\eqn{ef3}
$$
 q\le \frac{t\ef'(t)}{\ef(t)}\le p, \quad t>0.
$$
\end{prop}

If a Young function $\ef$ satisfies \eqref{ef3}, then the function
$t^{p}\ef(t)$ is nondecreasing and the function $t^q\ef(t)$ is nonincreasing. 
In this case the Orlicz space $L^{\ef}(\Omega)$ coincides with the set of all
measurable functions for which the integral $\int_{\Omega}\ef(|g|)\,dx$ 
converges.

In the following paragraphs we recall the definition of capacity
and certain pointwise properties of Sobolev functions.   We 
refer the reader to \cite{FZ} for the proofs.

For $1 \le p < \infty$ the $p$-capacity $\cp_p(E)$ of a set 
$E \subset \R^n$ is the infimum of numbers of the form 
\begin{equation} 
\label{def-cap} 
\int_{\R^n} (|u|^p + |\nabla u|^p) \, dx 
\end{equation} 
corresponding to all $u \in W^{1,p}(\R^n)$ with the property 
that $u \ge 1$ on a neighborhood of $E$.  
$p$-capacity is related to Hausdorff measure as follows:  if $p > 1$ 
and $\cp_p(E) = 0$, then $\hm^q(E) = 0$ for all $q > n-p$, whereas 
$\hm^{n-p}(E)< \infty$ implies $\cp_p(E) = 0$.  As for the case $p=1$, 
Fleming \cite{whf} proved that $\cp_1(E) = 0$ if and only if  
$\hm^{n-1}(E) = 0$.  Related statements concerning the Hausdorff
content are given in Lemma \ref{l:key-scalar} and
Theorem \ref{t:msz-cap} below.  

A function $f : \Omega \to \R^n$ is said to be $p$-quasicontinuous
if for every $\eps > 0$ there exists an open set $G$ with
$\cp_p(G) < \eps$ so that $f|_{\Omega \setminus G}$ is continuous.
Any two $p$-quasicontinuous functions which agree almost everywhere
agree up to a set $E$ with $\cp_p(E) = 0$.

If $\Omega\subset\R^n$ is an open set and 
$f\in L_{\loc}^{1}(\Omega)$, then a function
$\tilde f$ is said to be a \emph{precise representative}
of $f$ if 
\eqn{defpreciserep}
$$
	\tilde f(x) := \lim_{r \to 0} \, \intav{B(x,r)} f(y) \, dy
$$
at all points $x$ where this limit exists.  It is clear from the
Lebesgue differentiation theorem that
any function $f \in L^1_{\loc}(\Omega)$ may be modified on a set
of Lebesgue measure zero so as to be precisely represented and
that any two precise representatives coincide almost everywhere.
A mapping $f \in L^1_{\loc}(\Omega;\R^m)$ is said to be 
\emph{precisely represented} if each of its component functions
is a precise representative.

Fundamental properties of functions $f \in W^{1,p}(\Omega)$ are
that the limit (\ref{defpreciserep}) exists for all $x$ outside
a set $E$ with $\cp_p(E) = 0$, and that any precise representative
of a function $f \in W^{1,p}(\Omega)$ is $p$-quasicontinuous.
Thus if $p > m$ or $p = m = 1$,
any two precise representatives agree outside a set $E$ with
$\hm^{n-m}(E) = 0$. Since sets of $\hm^{n-m}$ measure zero
are negligible for results of Theorem \ref{t:main} and Theorem \ref{t:graph}, 
it follows that their statements
will hold for any $p$-quasicontinuous representative and even any
representative with Lebesgue points $\hm^{n-m}$ almost everywhere.

Similarly it is proven in \cite{MSZ} that
if $f \in W^{1,1}_{\loc}(\Omega)$ and $|\nabla f| \in L^{m,1}(\Omega)$,
then any precise representative of $f$ is defined up to a set $E$
with $\hm^{n-m}(E) = 0$.  
It follows that the statement of Theorem \ref{t:lorentz}
will hold for any
representative with Lebesgue points $\hm^{n-m}$ almost everywhere.

\section{Lusin's Condition (N)}

A mapping $f : \R^n \to \R^k$, $k \ge n$, is said to satisfy Lusin's
condition (N) if $\hm^n(f(E)) = 0$ whenever $E \subset \R^n$ satisfies
$\mathcal{L}^n(E) = 0$.  In this section we show that if 
$f \in W^{1,p}(\R^n;\R^m)$ and $\bar f : \R^n \to \R^{n+m}$ 
satisfies condition (N) then the conclusions of Theorems \ref{t:main}
and \ref{t:graph} hold.  This holds without any particular restriction
on $p$.  First we require two well known lemmas.

\begin{lemma}[Eilenberg inequality, \cite{eilenberg}]\label{eilenberg}
Suppose $m \le d \le m+n$, $A \subset \R^{n+m}$, and $h : A \to \R^m$ is
Lipschitz.  Then
\[
\int^{*}\hm^{d-m}(A\cap h^{-1}(y))\,d\hm^{m}(y)\leq C
(\lip h)^{m}\hm^{d}(A), 
\]
where $\int^{*}$ denotes the upper integral, $\lip f$ is the
Lipschitz constant of $f$, and $C=C(m)$ is a constant depending
only on $m$.
\end{lemma}

\begin{lemma} \label{l:count:lip}
If $f \in W^{1,1}_\loc(\R^n;\R^m)$ then there exist Lipschitz functions
$f_k\colon\R^n\to\R^m$ and disjoint subsets $E_k$ of $\R^n$ 
such that $f = f_k$ on $E_k$ and
$\mathcal{L}^n (\R^n \setminus \cup E_k) = 0$.
\end{lemma}

This follows from the a.e.\ approximate differentiability of $f$
\cite[3.1.4]{GMT}, \cite{BZ} and a general 
property of a.e.\ approximately differentiable functions
\cite[3.1.8]{GMT}.

\begin{theorem}
\label{t:cond(N)}
Let $f \in W^{1,1}_{\loc}(\R^n;\R^m)$, $1 \le m \le n$, and
suppose that $\bar f$ satisfies condition (N). 
Then $f^{-1}(y)$ is countably $\mathcal{H}^{n-m}$ rectifiable
for almost all $y \in \R^m$,
the graph $\G_f$ is countably $\mathcal{H}^n$ rectifiable and 
the coarea formula \eqref{coarea2} and area formula \eqref{area2}
hold for all measurable sets $E \subset \Omega$. 
\end{theorem}

\begin{proof}  Choose functions $f_k$ and sets $E_k$ 
as in Lemma~\ref{l:count:lip}.  Since $f_k=f$ on $E_k$ and
$\nabla f_k=\nabla f$ a.e.\ on $E_k$, it follows from the
classical area
 formula for Lipschitz functions 
\cite[3.2.5]{GMT} that that
\begin{equation}
\label{e:area}
\hm^n(\bar f(E)) = \int_E |J_n \bar f(x)| \, dx
\end{equation}
whenever $E \subset E_k$.  Moreover (\ref{e:area}) will be satisfied
whenever $\mathcal{L}^n(E) = 0$ provided that $\bar f$ satisfies condition
(N).  Therefore the monotone convergence theorem implies (\ref{e:area})
for all $E \subset \R^n$, as desired.  
That the graph $\mathcal{G}_f$ is countably $\hm^n$ rectifiable 
is evident since
\[
\hm^n\left( \mathcal{G}_f \setminus \bigcup_{k=1}^\infty
\bar f(E_k) \right) = 0.
\]

Now denote by $\pi : \R^{n+m} \to \R^n$ and $\rho : \R^{n +m} \to \R^m$ the
projections
\[
\pi(x, x_{n+1}, \ldots, x_{n+m}) = x \text{ and }
\rho(x, x_{n+1}, \ldots, x_{n+m}) = (x_{n+1}, \ldots, x_{n+m}).
\]
Let $E \subset \R^n$ satisfy $\mathcal{L}^n(E) = 0$ and apply the
Eilenberg inequality with $A = \bar f(E)$, $h = \rho$, and
$d=n$.  Then
\[
\int_{\R^m}\hm^{n-m}(A \cap \rho^{-1}(y)) \,d \hm^m(y) \le
C \hm^n(\bar f(E)) = 0.
\]
Now, $A \cap \rho^{-1}(y) = \{(x,f(x)) \in \R^{n+m} 
\colon x \in E, \, f(x) = y\}$, so
it follows that $\pi(A \cap \rho^{-1}(y)) = E \cap f^{-1}(y)$.  Since
Hausdorff measures do not increase on projection we conclude
\begin{equation}
\label{co-N}
\int_{\R^m} \hm^{n-m}(E \cap f^{-1}(y)) \, d \hm^m(y) = 0.
\end{equation}
From (\ref{coarea2}) above we have that
\begin{equation}
\label{e:coarea}
\int_{\R^m} \hm^{n-m}(E \cap f^{-1}(y)) \, d \hm^m (y) = 
\int_E |J_m f(x)| \, dx
\end{equation}
holds whenever $E \subset E_k$, so again the monotone convergence 
theorem implies that (\ref{e:area}) holds for all $E \subset \R^n$.
That $f^{-1}(y)$ is countably $\hm^{n-m}$ rectifiable for almost
$y \in \R^m$ follows from the fact that for each $k$, 
$E_k\cap f^{-1}(y)=E_k\cap f_k^{-1}(y)$ is countably $\hm^{n-m}$ 
rectifiable for almost all $y \in \R^m$ (see~\cite[3.2.15]{GMT}), 
and that $\hm^{n-m}(E \cap f^{-1}(y)) = 0$ for
almost all $y$ whenever $\mathcal{L}^n(E) = 0$.  
\end{proof}

\section{A general criterion for condition $N$}

Throughout this section we denote by
$\pi : \R^{n+m} \to \R^n$
the projection
$$
\pi(x, x_{n+1}, \ldots, x_{n+m}) = x.
$$

\begin{lemma}
\label{l:lemma5} Suppose $m\le d\le m+n$.
Let $E \subset \R^{n+m}$.  Then 
$$
\hm^{d}_{\infty}(E) \le C (\diam \, E)^m
\hm^{d-m}_\infty (\pi(E)).
$$
where $C = C(m,n,d)$.
\end{lemma}

\begin{proof}
Assume $\diam \, E < \infty$. 
Cover $\pi(E)$ with a sequence of sets
$\{A_k\}_{k=1}^\infty$ of finite and positive diameter.
For each $k$ let $N_k$ denote the greatest
integer satisfying 
\[
(N_k - 1) \cdot (\diam \, A_k) < \diam \, E.
\]
The set $E \cap \pi^{-1}(A_k)$ is easily seen to be contained
in a cylinder of the form $A_k \times Q_k$, where $Q_k$ is
a cube in $\R^m$ with sidelength $\diam \, E$.  Since
$N_k \cdot (\diam \, A_k) \ge \diam \, E$, $Q_k$ may be
covered by $N_k^m$ cubes $\{Q_k^j\}$ of sidelength 
$\diam \, A_k$.  It follows that
\begin{eqnarray*}
\hm^{d}_\infty(E \cap \pi^{-1}(A_k)) & \le & C \sum_{j=1}^{N_k^m}
\left( \frac {\diam\,(A_k \times Q_k^j)}{2} \right)^{d} \\
& \le & C N_k^m \cdot (\diam\, A_k)^{d} \\
& \le & C (\diam \, A_k + \diam\, E)^m \cdot
(\diam A_k)^{d-m}.
\end{eqnarray*}
Since $\diam \, A_k \le \diam \, \pi(E) \le \diam \, E$, summing
over $k$ shows that
\[
\hm^{d}_\infty(E) \le C (\diam \, E)^m \sum_{k=1}^\infty
(\diam \, A_k)^{d-m}.
\]
The result follows by taking the infimum over all coverings $\{A_k\}$.
\end{proof}

\begin{cor}\label{c:repres}
Suppose $m\le d\le m+n$.
Let $E \subset \R^{n+m}$, $\hm^{d-m}(\pi(E))=0$. Then
$\hm^{d}(E)=0$. 
\end{cor}

The following lemma concerns a measurability property of the
graph of a mapping $f$ and is used in the proof of Theorem
\ref{t:general} below.  Recall that $\bar f$ denotes the graph
mapping $\bar f(x) = (x,f(x))$ and $\mathcal G_f$ the graph set
$\bar f(\R^n) \subset \R^{n+m}$.

\begin{lemma}
\label{l:measure}
Let $f : \R^n \to \R^m$ be an $\L^n$ measurable mapping.  Then
$\pi (\mathcal{G}_f \cap E)$ is $\L^n$ measurable for every
Borel set $E \subset \R^{n+m}$.
\end{lemma}

\begin{proof}
Let $f^*$ be a Borel measurable representative of $f$.
The graph $\mathcal{G}_{f^*}$ of such a mapping is a Borel subset of
$\R^{n+m}$, implying that the projection $\pi(\mathcal{G}_{f^*} \cap E)$ 
is a Borel subset of $\R^n$ for any Borel set $E \subset \R^{n+m}$.
See \cite[\S 31.\ VII]{K}.  Since $f$ and $f^*$ agree up to a set
of $\L^n$ measure zero, so do the sets $\pi(\mathcal{G}_{f^*} \cap E)$ 
and $\pi(\mathcal{G}_f \cap E)$.  Thus $\pi(\mathcal{G}_f \cap E)$
is $\L^n$ measurable.
\end{proof}

The following result yields a criterion similar to that of Rad\'o and
Reichelderfer, see \cite{RR}, \cite{M-ac}.

\begin{theorem}\label{t:general} Suppose that $1\le m\le n$.
Let $f \in W^{1,1}_{\loc}(\R^n;\R^m)$ and suppose that there is
$\theta \in L^1_{\loc}(\R^n)$ such that 
\begin{equation}
\label{general}
      \hm^{n-m}_\infty \bigl( \pi(\mathcal{G}_f \cap B(z,r)) \bigr) \le r^{-m}
      \int_{\pi(\mathcal{G}_f \cap B(z,4r))} \theta(x) \, dx
\end{equation}
for all $z \in \R^{n+m}$ and $r > 0$.  Then there is a constant
$C = C(m,n)$ depending only on $m$ and $n$ so that
\[
\hm^n (\bar f(E)) \le C \int_E \theta(x) \, dx.
\]
for all $\L^n$ measurable $E \subset \R^n$.
In particular, $\bar f$ satisfies condition (N).
\end{theorem}

\begin{proof}
Define a set function $\sigma$ on $\R^{n+m}$ by
\[
\sigma(E) = \int_{\pi(\mathcal{G}_f \cap E)} \theta (x) \, dx.
\]
Lemma \ref{l:lemma5} with $d=n$ and (\ref{general}) imply that
\begin{align*}
\hm^n_\infty \bigl( \mathcal{G}_f \cap B(z,r) \bigr) &\le
C  r^m\hm^{n-m}_\infty \bigl( \pi( \mathcal{G}_f \cap B(z,r))
\bigr) \\
&\le C  \int_{\pi (\mathcal{G}_f \cap B(z,4r))} \theta(x) \, dx
\\ 
&=  C  \sigma (B(z,4r))
\end{align*} 
for any $z \in \R^{n+m}$ and $r > 0$.
  Using the fact that 
\[
\limsup_{r \to 0^+} r^{-n} \hm^n_\infty(\mathcal{G}_f \cap B(z,r)) \ge C
\]
for $\hm^n$ almost every $z \in \mathcal{G}_f$
(\cite[Lemma 10.1]{F5}),  it follows that
\begin{equation}
\label{lower-est}
\limsup_{r \to 0^+} r^{-n} \sigma(B(z,r)) \ge C 
\end{equation}
for $\hm^n$ almost all $z \in \mathcal{G}_f$.  Now, Lemma \ref{l:measure}
implies that $\sigma$ is a measure on the Borel $\Sigma$-algebra of
$\R^{n+m}$, and so it may be extended to a regular Borel (outer) measure 
$\sigma^*$ on all of $\R^{n+m}$ in the usual way:
\[
	\sigma^*(E) := \inf \left\{ \sigma(B) : \, E \subset B, \,
	\text{$B$ a Borel set} \right\}.
\]
Since $\theta$ is locally integrable it follows that $\sigma^*$ is
a Radon measure on $\R^{n+m}$.  Therefore~
\cite[Theorem 10.3]{F5},
\cite[Theorem 6.9]{M}, and (\ref{lower-est}) 
imply that
 $\hm^n(E) \le C \sigma^*(E)$ for all
$E \subset \mathcal{G}_f$
. 
 Finally, given an $\L^n$ measurable set
$E \subset \R^n$, choose a Borel
set $G$ with $E \subset G$.  Then $\bar f(E) \subset G \times \R^m$,
$G \times \R^m$ is Borel, and
\[
\hm^n(\f(E)) \le C \sigma^*(\bar f (E))
\le C \sigma(G \times \R^m) 
= C \int_G \theta(x) \, dx
.
\]
The proof is completed by taking the infimum over all such $G$.
\end{proof}

\section{Capacitary estimates}

The next lemma gives a capacitary estimate for $W^{1,1}$ which
relates the $1$-capacity with the $(n-1)$-dimensional Hausdorff
content.  This result is due to Fleming~\cite{whf},
whose proof depends on information based on BV
functions, in particular sets of finite perimeter. We give an
independent proof which circumvents
the need for BV theory. 

\begin{lemma}\label{l:key-scalar}
Let $E\subset \R^n$. Then $\hm_{\infty}^{n-1}(E)\le C \cp_{1} (E)$,
where $C = C(n)$.
\end{lemma}

\begin{proof}
Let $u$ be a test function for the capacity.
Let us consider the continuous monotone real function
$$
\psi(t) = \int_{\{0<u<t\}}(|u|+ |\nabla u|)\,dy.
$$
Then $\psi$ is a.e.\ differentiable and
$$
\int_0^1\psi'(t)\,dt\le \psi(1)-\psi(0)=\psi(1).
$$
Hence there exists $s\in (0,\ell)$ such that
$\psi'(s) < 2\psi(1)$.
We find $\delta>0$ such that
\eqn{howdelta}
$$
\frac{\psi(s')-\psi(s)}{s'-s}\le 2\psi(1)
\quad\text{for each }s'\in (s,s+\delta).
$$
Let $x\in E$. Then $u(x)\ge 1$ in a neighborhood of $x$ and so
the set
$$
R_x:=\{\rho:\; |B(x,\rho)\cap \{u\ge s\}| < \tfrac12 |B(x,\rho)|\}
$$
is nonempty. Let us define
$$
r_x = \sup R_x.
$$
Then obviously $\sup_{x\in E}r_x<\infty$ and
\eqn{med-upper}
$$
|B(x,r_x)\cap\{u<s\}|\le \frac12 |B(x,r_x)|.
$$
and
\eqn{med-lower}
$$
|B(x,r_x)\cap\{u\ge s\}|\le \frac12 |B(x,r_x)|.
$$
We use the Besicovitch covering theorem to extract a (finite or infinite)
sequence $\{B_j\}_{j\in I}$ of balls $B_j=B(x_j,r_j)$ from $\{B(x,r_x)$:
$x\in E\}$
such that it covers the set
$E$ and its overlap
multiplicity is bounded by an integer $M$ depending only on $n$.
Here $I=\en$ or $I=\{1,2,\dots,i_{\max}\}$.
Fix $i\in I$. Using \eqref{med-lower},
we find a level $t_i\in (s,s+\delta)$ such that
\eqn{med-lower3}
$$
|B(x_j,r_j)\cap\{u\le t_i\}|\le \tfrac 34 |B(x_j,r_j)|, \quad j=1,\dots,i.
$$
Consider the truncated function
$$
u_i =
\begin{cases}
t_i& \text{ on }\{u\ge t_i\},\\
s &\text{ on }\{u\le s\},\\
u &\text{ on }\{s< u<t_i\}.
\end{cases}
$$
Let $s_{i,j}$ be a median value of $u_i$ on $B_j$. From
\eqref{med-upper} and \eqref{med-lower3} we infer that
$$
t_i-s\le C\vint_{B_j}|u_i-s_{i,j}|\,dy.
$$
By the Poincar\'e inequality we can continue
$$
t_i-s\le Cr_j\vint_{B_j}|\nabla u_i|\,dy.
$$
This means
$$
r_j^{n-1}\le \frac{C}{t_i-s}\;\int_{B_j}(|u_i|+|\nabla u_i|)\,dy.
$$
Summing over $j=1,\dots,i$ and using \eqref{howdelta} we obtain
$$
\aligned
\sum_{j=1}^ir_j^{n-1}&\le
\frac{C}{t_i-s}\;\sum_{j=1}^i\int_{B_j}(|u_i|+|\nabla u_i|)\,dy
\le \frac{CM}{t_i-s}\int_{B}(|u_i|+|\nabla u_i|)\,dy
\\&
= \frac{CM}{t_i-s} \int_{B\cap \{s\le u<t_i\} } 
(|u| + |\nabla u|) \, dy = CM
\frac{\psi(t_i)-\psi(s)}{t_i-s}
\\&\le C\psi(1)M.
\endaligned
$$
We can pass with $i$ to $i_{\max}$ or $\infty$ and we obtain
$$
\hm_{\infty}^{n-1}(E)\le C \psi(1)\le \int_{\R^n}(|u|+|\nabla u|)\,dy
$$
as required.
\end{proof}

\begin{cor}
\label{c:scalar-cap}
Let $E\subset\R^n$. Suppose that $u\in W^{1,1}(\R^n)$
 is
precisely represented and that $u\ge 1$ on $E$.
  Then
$$
\hm_{\infty}^{n-1}(E)\le C\int_{\R^n}\bigl(|\nabla u|+ |u|\bigr)\,dx.
$$
\end{cor}

\begin{proof} Any precise representative of 
$u\in W^{1,1}(\R^n)$ has
Lebesgue points $\hm^{n-1}$-a.e., cf.\ \cite{whf}, \cite{FZ}.
Hence we may assume that $u$ has a Lebesgue point at 
$\hm^{n-1}$-a.e.\ point of $x \in E$.  Given $\eps>0$,
let $\{u_k\}$ be a sequence of smooth function approximating $u$ and 
obtained by mollification, such that 
$$
\|u_k-u\|_{1,1}\le 2^{-k-1}\eps.
$$
Set 
$$
w = u_1+\sum_{k=1}^{\infty}|u_{k+1}-u_k|.
$$
Then $w\in W^{1,1}(\R^n)$, $\|w\|_{1,1}\le \|u\|_{1,1}+\eps$ and
$w \ge 1$ on a neighborhood of $E$. Hence
$$
\hm_{\infty}^{n-1}(E)\le 
C\int_{\R^n}\bigl(|\nabla w|+ |w|\bigr)\,dx
\le C\int_{\R^n}\bigl(|\nabla u|+ |u|\bigr)\,dx + C\eps.
$$
Letting $\eps\to 0$ we obtain the assertion.
\end{proof}

For the proofs of Theorems \ref{t:main} and \ref{t:graph} we need the 
well known relation between the $p$-capacity and the Hausdorff content,
see e.g.\ \cite[Theorem 5.1.13]{AH}. 
This result goes back to Frostman ($p=2$) and the case of general $p$
is due to 
Reshetnyak \cite{R} and Maz'ya and Havin \cite{MH}.

\begin{theorem}
\label{t:haus-cap}
Suppose that $1\le m< p$ and $E \subset \R^n$.  Then
$$
\hm_{\infty}^{n-m}(E)\le C \cp_p(E),
$$
where $C = C(n,m,p)$.
\end{theorem}

Since precise representatives of functions $u \in W^{1,p}(\R^n)$ have
Lebesgue points $\hm^{n-m}$ almost everywhere \cite{FZ} we may continue
as in the proof of Corollary \ref{c:scalar-cap} to obtain the following.

\begin{cor}
\label{c:haus-cap}
Suppose that $1\le m< p$, $E\subset\R^n$, $u\in W^{1,p}(\R^n)$
is precisely represented, and that $u\ge 1$ on $E$.  Then
$$
\hm_{\infty}^{n-m}(E)\le C \int_{\R^n}\bigl(|\nabla u|^p+ |u|^p\bigr)\,dx.
$$
\end{cor}

Finally, to reach the generality of Theorem \ref{t:lorentz} we need the
following capacitary estimate derived in \cite{MSZ}.  

\begin{theorem}\label{t:msz-cap}
Suppose that $m>1$ is an integer.
Let $\ef$ be a Young function satisfying \eqref{ef2} and \eqref{ef3}.
Let $E\subset \R^n$.
Suppose that $u\in W_{\loc}^{1,1}(\R^n)$
is precisely represented, $u\ge 1$ on $E$.
Then
$$
\hm_{\infty}^{n-m}(E)\le C\int_{\R^n}\bigl(\ef(|\nabla u|)+\ef(|u|)\bigr)\,dx.
$$
\end{theorem}

\section{Application of the general criterion}

In this section we prove Theorem \ref{t:lorentz}. The reader interested only
in the Lebesgue scale of spaces may read the arguments with $\ef(t)=t^p$
and appeal to Corollary \ref{c:haus-cap}.  This is sufficient to establish
Theorems \ref{t:main} and \ref{t:graph}.  To prove the more general result we
use the estimate in Theorem \ref{t:msz-cap}. 

\begin{proof}
In view of Theorem \ref{t:cond(N)} it is sufficient to verify the
assumptions of Theorem \ref{t:general}.
Select a point $z \in \R^{n+m}$ and $r > 0$.  Writing $z = (x_0, y_0)$
we have
\[
\mathcal{G}_f \cap B(z,r) \subset \mathcal{G}_f \cap
\left[ B(x_0,r) \times B(y_0,r) \right]
\]
hence
\[
\pi \left( \mathcal{G}_f \cap B(z,r) \right) \subset
B(x_0,r) \cap f^{-1}(B(y_0,r)).
\]
Let $E = B(x_0,r) \cap f^{-1}(B(y_0,r))$ and define $\hat E$ by
\[
\hat E = \frac{1}{2r} (E - x_0) = \{x \in \R^n : x_0 + 2rx \in E\}.
\]
Then
\[
\hm^{n-m}_\infty(\hat E) = (2r)^{m-n} \hm^{n-m}_\infty (E).
\]
and thus
\begin{equation}
\label{scaling}
\hm^{n-m}_\infty \bigl(\pi (\mathcal{G}_f \cap B(z,r))\bigr) \le 
C r^{n-m} \hm^{n-m}_\infty(\hat E),
\end{equation}
where $C = C(m,n)$.  
Now, $\xi \in \hat E$ implies that
\[
|\xi| \le \frac{1}{2} \quad\text{and}\quad
\frac{|f(x_0 + 2r\xi) - f(x_0)|}{2r} \le \frac{1}{2}.
\]
Thus we consider the test function $u\eta$ where
$$ 
u(\xi) = 2 \Bigl(1- 
\frac{|f\bigl(x_0+ 2r\xi)\bigr)-f(x_0)|}{2r}\Bigr)^+ 
$$ 
and $\eta$ is a smooth cutoff function such that $\chi_{\strut 
B(0,\tfrac12)}\le \eta\le \chi_{\strut B(0,1)}$. 
If $m=1$, we set $\ef (t) = t$ and apply Corollary \ref{c:scalar-cap}
to the function $u\eta$.  If $m > 1$, we use Proposition
\ref{p:msz} to find a Young function $\ef$ satisfying (\ref{ef2})
and (\ref{ef3}) such that
$$
\int_\Omega \ef (|\nabla u|) \, dx < \infty
$$
and apply Theorem \ref{t:msz-cap} to the function $u\eta$.  In 
either case we obtain
$$
\aligned
\hm^{n-m}(\hat E) &\le  C \int_{B(0,1)}
\bigl(\ef(| u\eta |)+ \ef(|\nabla( u\eta ))\bigr)\,d\xi
\\&
\le C \int_{B(0,1)\cap \{u>0\}}
\bigl(1+ \ef(|\nabla u|)\bigr)\,d\xi
\endaligned
$$
Applying the change of variable $x \to x_0 + 2r\xi$ this becomes
\begin{equation}
\label{capacit}
\hm^{n-m}(\hat E) \le C r^{-n} \int_{B(x_0, 2r) \cap f^{-1}(B(y_0,2r))}
\bigl(1+ \ef(|\nabla u|)\bigr) \, dx.
\end{equation}
Since 
\[
B(x_0, 2r) \cap f^{-1}(B(y_0,2r)) \subset \pi(\mathcal{G}_f
\cap B(z,4r))
\]
(\ref{scaling}) and (\ref{capacit}) above imply
$$
\hm^{n-m}_\infty\left( \pi(\mathcal{G}_f \cap B(z,r)) \right)
\le C \,r^{-m}\int_{\pi(\mathcal{G}_f \cap B(z,4r))}
\bigl(1+ \ef(|\nabla u|)\bigr)\,dx,
$$
verifying the assumptions of Theorem \ref{t:general} with 
$\theta = C\bigl(1+ \ef(|\nabla u|)\bigr)$ and thus concluding the proof.
\end{proof}

We now establish a result that yields Theorem \ref{t:main} under a
condition that allows some of
the coordinate functions to be members of $W^{1,p}(\R^n)$, $p<m$,
provided that the remaining ones are Lipschitz. For this purpose,
let $1\leq k<m<n$ and let $h\colon\R^n\to\R^{m}$ have the form
$$
	h = (f_{1},\ldots,f_{k},g_{k+1},\ldots, g_{m})
$$
with $f_{1},\ldots,
f_{k}\in W^{1,p}(\R^n), p>k$, and $g_{k+1},\ldots,g_{m}$ Lipschitz.

\begin{theorem} With $h$ as above, $|Jh|\in L^{1}(\R^n)$ and 
$$
\int_{\R^n}|Jh|\;dx=\int_{\R^{m}}\hm^{n-m}(h^{-1}(y))\,dy
$$
for each measurable set $E \subset \R^n$.
\end{theorem}

\begin{proof}
By Theorem \ref{t:cond(N)} it suffices to show that $\hm^n(\bar f(N)) = 0$
whenever $N \subset \R^n$ and $\mathcal{L}^n(N) = 0$.  Write $h = (f,g)$
with $f = (f_1, \ldots, f_k)$ and $g = (g_{k+1}, \ldots, g_m)$ and let
$N \subset \R^n$ satisfy $\mathcal{L}^n(N) = 0$.  Identify $\R^{n+k}$
with $\R^n \times \R^k$ and write $z \in \R^{n+k}$ as $z = (z',z'')$ with
$z' \in \R^n$, $z'' \in \R^k$.  By Theorem \ref{t:graph} the set
$\bar f(N) \subset \R^{n+k}$ has $\hm^n$ measure zero.  Define
$g^* : \R^{n+k} \to \R^{n+m}$ by $g^*(z) = (z,g(z'))$.  Then $g^*$ is
Lipschitz, $h = g^* \circ \bar f$, and therefore $h(N) = g^*\circ \bar f(N)$
has $\hm^n$ measure zero in $\R^{n+m}$.
\end{proof}

\section{H\"older continuous mappings}

In this section we prove Theorem~\ref{t:hoelder}.  The following
result from \cite{HMy} plays a crucial role in the proof.

\begin{theorem}\label{p:hmy}
Let $p<n$, $a>0$, $\beta\in(0,1)$ and $\gamma>0$.
Let $\Omega\subset\R^n$ be an open set
and $u\in W^{1,m}(\Omega)$ be a nonnegative function 
such that $u>0$ a.e.
Let 
\eqn{zz}
$$
Z=\Bigl\{z\in \Omega: 
\limsup_{r\to 0}r^{-\beta}\intav{B(z,r)}u\,dx<\gamma\Bigl\}.
$$
Suppose that $\hm^{n-m}(Z)>a$. 
Then there exists a compact set $F\subset \Omega\setminus Z$
such that $\L^n(F)>0$ and
$$
\sup_{F}u^m\le C\int_{F}|\nabla u|^m\,dx,
$$
where $C=C(n,m, a,\beta)$.
\end{theorem}

\begin{proof}[Proof of Theorem~\ref{t:hoelder}]
As in the proof of Theorem~\ref{t:main} above it is sufficient
to verify (\ref{co-N}) whenever $E \subset \R^n$ and
$\L^n(E) = 0$.  Find an open set $\Omega$ such that $E \subset \Omega$
and
\eqn{intsmall}
$$
\int_{\Omega}|\nabla f|^m\,dx<\varepsilon .
$$
Set
$$
Y=Y_{a}=\bigl\{y\in \erm: \
\L^n(E\cap f^{-1}(y))=0,\
\hm^{n-m}(E\cap f^{-1}(y))>a
\bigr\}.
$$
Owing to the H\"{o}lder continuity of $f$,
given $y\in Y$ we may apply Proposition \ref{p:hmy}
to the set
$\Omega$ and the function $u=|f-y|$ and find a $\delta(y)>0$ 
and a compact set $F(y)\subset\Omega$ such that $\L^n(F(y))>0$ and
\eqn{applhmy}
$$
\sup_{F(y)}|f(x)-y|^m <\delta(y) \le \int_{F(y)}|\nabla f|^m\,dx.
$$
Using a Vitali type covering argument to the system of balls
$B(y,\delta(y))$ we find disjointed 
balls $B(y_j,\delta_j)$ and sets $F_j$
such that $y_j\in Y$, $\delta_j=\delta(y_j)$, $F_j=F(y_j)$,
and
$$
Y\subset\bigcup_j B(y_j,5\delta_j).
$$
Since the balls $B_j$ are disjointed, the sets $F_j\subset f^{-1}
(B(y_j,\delta_j))$ are also disjointed. We infer that
$$
\aligned
\L^m(Y)&\le \sum_j\L^m(B(y_j,5\delta_j))
\le C\sum_j\delta_j^m\le C\sum_j\int_{F_j}|\nabla f|^m\,dx
\\&
\le C\int_{\Omega}|\nabla f|^m\,dx<C\varepsilon .
\endaligned
$$
Letting $\varepsilon \to0$ we obtain that $\L^m(Y)=0$. Since $a>0$ was
arbitrary and the set
$$
\bigl\{y\in \erm:\
\L^n(f^{-1}(y))>0   
\bigr\}
$$
is countable (and thus of zero $m$-dimensional measure), we easily 
conclude the proof.
\end{proof}

\begin{rem}
\label{r:proof}
The proof above can be modified to give an independent proof
of Theorem~\ref{t:main} when $p > m$.  Indeed, by \cite[ Prop. 3.2]{HMy}. 
$\hm^{n-m}$-a.e.\ $z \in \R^n$ satisfies
\[
\limsup_{r \to 0} r^{-\beta} \intav{B(z,r)}|f(x) - f(z)| \, dx <
\infty 
\]
with $\beta = 1 - \frac{m}{p}$, provided that $f \in W^{1,p}(\R^n;\R^m)$ and
$p > m$.

Alternatively, when $p > 1$ one can use the fact that for $0 < \lambda < 1$
and $\eps > 0$ there exists an open set $U$ and a mapping 
$g \in W^{1,p}(\R^n;\R^m)$ so that $g$ is H\"{o}lder continuous with
exponent $\lambda$, $f(x) = g(x)$ for all $x \in \R^n \setminus U$, 
and $B_{1-\lambda,p}(U) < \eps$, where $B_{1-\lambda,p}$ is the Bessel
capacity.  See \cite{BHS}, \cite{S}; a weaker but also sufficient
result of this type is given in \cite{m-holder}.
Thus for $f \in W^{1,p}(\R^n;\R^m)$ there is a set $N$ with
$B_{1-\lambda,p}(N) = 0$ so that (\ref{coarea2}) holds for all 
$E \subset \R^n \setminus N$.  In case $p > m$ then $\lambda$ may
be chosen so that $\hm^n(N) = 0$, extending (\ref{coarea2}) to all 
measurable sets $E \subset \R^n$.
\end{rem}

\section{Further results on the graph mapping}

In this section we consider properties of the graph $\mathcal{G}_f$
for mappings in the borderline case $W^{1,m}(\R^n;\R^m)$, for
$1 \le m \le n$.  The following Theorem was proved in~\cite{MM} for 
$k = n$ and extended to $k \ge n$ in~\cite{M-ag}.

\begin{theorem}
\label{t:area}
Suppose that $k \ge n$ and that $f \in W^{1,n}_{\loc}(\R^n;\R^k)$ is
H\"{o}lder continuous.  Then $f$ satisfies condition (N).
\end{theorem}

We observe that mappings $f \in W^{1,p}(\R^n;\R^k)$, $p > n$, satisfy the
hypothesis of Theorem~\ref{t:cond(N)} and therefore this result gives
a new proof of the area formula of Marcus and Mizel~\cite{MMi}.

We let $\grass(n,m)$ denote the Grassmann manifold of $m$ dimensional 
subspaces of $\R^n$.  Given $V \in \grass(n,m)$ let $V^\perp \in 
G(n,n-m)$ denote its orthogonal complement.  Writing $x \in \R^n$ 
as $x = x_V + x_{V^\perp}$, with $x_V \in V$ and $x_{V^\perp} \in 
V^\perp$, the orthogonal projection $P_V : \R^n \to V$ is given by 
$P_V(x) = x_V$.  If $V \in \grass(n,m)$ and $W \in \grass(n,k)$, 
$m+k \ge n$, observe that 
\[ 
V^\perp \subset W \text{ if and only if } 
W^\perp \subset V. 
\] 
In this case we write $V \perp W$.  The following proposition is 
an immediate consequence of the definitions. 
 
\begin{prop}\label{p:ortho} 
If $V \perp W$ then  
\begin{enumerate} 
\item $P_V(W) = V \cap W$ 
\item $P_V(W + a) \cap  P_V(W + b) = \emptyset$ for all $a,b \in W^\perp$,  
$a \not= b$. 
\end{enumerate} 
\end{prop} 

Recall that the $n$ dimensional integralgeometric measure of a set $E \subset
\R^{n+m}$ is given by
\[
\mathcal{I}^n(E) = \int_{\grass(n+m,n)} \int_V N(P_V,E,y)\,
d\hm^n(y) \, d\sigma_{n+m,n}(V),
\]
where $N(P_V,E,y)$ is the number of points $x \in E$ satisfying 
$P_V(x) = y$, and $\sigma_{n+m,n}$ is a normalized Haar measure on
$\grass(n+m,n)$.  In order to prove Theorem~\ref{t:hoelderarea} it suffices to
show that $\mathcal{I}^n( \bar f(E)) = 0$ whenever $\L^n(E) = 0$.  This is a
consequence of the following theorem.  

\begin{theorem}\label{l:slices} 
Assume that $1 \le m \le n$ and that $f \in W^{1,m}(\R^n;\R^m)$ is
H\"{o}lder continuous.  If $N \subset \R^n$ and $\L^n(N) = 0$, 
then $\hm^n(P_V \circ \bar f(N)) = 0$ for every 
$V \in \grass(n+m,n)$. 
\end{theorem} 
 
\begin{proof} 
Identify $\R^n$ with the subspace $\{x \in \R^{n+m} : x_{n+1} = \ldots = 
x_{n+m}=0\}$ of $\R^{n+m}$ and choose $V \in \grass(n+m,n)$.  Let $W \in 
\grass(n+m,2m)$ be any subspace containing $(\R^n)^\perp$ and $V^\perp$. 
Then $\R^n \perp W$, hence $W^\perp \subset \R^n$.  Writing 
\[ 
\R^n = (\R^n \cap W) \oplus (\R^n \cap W^\perp) =  
(\R^n \cap W) \oplus W^\perp, 
\] 
it follows from Proposition~\ref{p:ortho} that the sets $\{ \R^n \cap (W + t)$: 
$t \in W^\perp\}$ are disjoint $m$ planes in $\R^n$ whose union is $\R^n$. 
Similarly, the sets $\{V \cap (W+t)$: $t \in W^\perp\}$ are disjoint $m$ planes 
in $V$ whose union is $V$.  
 
For brevity we write $W_t = W + t$.  Define $h : V \to W^\perp$ by  
$h = P_{W^\perp}|_V$.  Then for $t \in W^\perp$ we have $h^{-1}(t) =  
V \cap W_t$.  By the coarea formula for Lipschitz mappings between 
rectifiable sets~\cite[3.2.22]{GMT} it follows that  
\[ 
\int_E |J_n h| \, d\hm^n(x) = \int_{W^\perp} \hm^m (E \cap W_t) \, 
d\hm^{n-m}(t)  
\] 
for any $\hm^n$ measurable subset $E$ of $V$, 
where $J_n h$ is the (non-zero constant) $n$ dimensional Jacobian of $h$. 
Denoting this constant by $J$ and taking 
$E = P_V \circ \bar{f}(A)$ for a Borel set $A \subset \R^n$ we have 
\begin{align*} 
\hm^n(E) & = J^{-1} \int_{W^\perp} \hm^m(( P_V \circ 
\bar{f}(A)) \cap W_t) \, d\hm^{n-m}(t) \\ 
& = J^{-1} \int_{W^\perp} \hm^m((P_V(\bar{f}(A) \cap W_t)) 
\, d\hm^{n-m}(t) \\ 
& \le J^{-1} \int_{W^\perp} \hm^m(\bar{f}(A) \cap W_t)) \, 
d\hm^{n-m}(t) 
\end{align*}   
by Proposition~\ref{p:ortho} and the fact that Hausdorff measure does 
not increase under orthogonal projection. 
 
On the other hand,  
since $\bar{f} \in W^{1,m}_{\loc}(\R^n;\R^{n+m})$ and $W^{1,m}$ is 
invariant under any nonsingular linear change of coordinates, Fubini's theorem 
implies that 
\begin{equation} 
\label{f:slice} 
\bar{f}_t := \bar{f}|_{\R^n \cap W_t} \in  
W^{1,m}_{\loc}(\R^n \cap W_t ;\R^{n+m}) 
\end{equation} 
for $\hm^{n-m}$ almost all $t \in W^\perp$.  Since $m = \dim (\R^n \cap 
W_t)$,  Theorem~\ref{t:area} above implies that 
\[ 
\hm^m(\bar{f}_t(A)) = \int_{A \cap W_t} |J_m \bar{f}_t(x)| \, 
d\hm^m(x) 
\] 
for $\hm^{n-m}$ almost all $t \in W^\perp$ whenever $A \subset \R^n$ is Lebesgue 
measurable and where $J_m \bar{f}_t$ is the $m$ dimensional Jacobian of 
$\bar{f}_t$.  Proposition~\ref{p:ortho} implies that 
$\bar{f}_t(E) = \bar{f}(E) \cap W_t$, hence for all such $t$ we have 
\[ 
\hm^m(\bar{f}(A) \cap W_t) = \int_{A \cap W_t} |J_m 
\bar{f}_t(x)|  \, d\hm^m(x). 
\] 
 
To complete the argument let $N \subset \R^n$ be a set with Lebesgue measure 
zero and let $A \subset \R^n$ be a Borel set containing $N$ with Lebesgue 
measure zero.  Fubini's theorem implies that $\hm^m(A \cap W_t) = 0$ for 
$\hm^{n-m}$ almost all $t \in W^\perp$, and therefore 
\[ 
\hm^m(\bar{f}(A) \cap W_t) = 0 
\] 
for all such $t$.  It follows that 
\[ 
\hm^n(P_V \circ \bar{f}(N)) \le \hm^n(P_V \circ \bar{f}(A)) 
\le J^{-1} \int_{W^\perp} \hm^m (\bar{f}(A) \cap W_t) \, d\hm^{n-m}(t) 
= 0, 
\] 
establishing the lemma.   
\end{proof} 
 
This result is not sufficient to conclude that $\hm^{n}(\f(N))=0$, 
since there exist sets of positive Hausdorff measure whose projection
in every direction has zero measure, cf. \cite{M}[Example 9.2]. 
However, with the help of \cite[Corollary 9.8]{M} and 
\cite[3.2.27]{GMT} we arrive at the following.
\begin{cor} With $f$ and $N$ as in the proof of Theorem \ref{l:slices},
$\mathcal{I}^n (\f(N))=0$. In
particular, the Hausdorff dimension of $\f(N)$ does not exceed $n$.  Moreover,
either $\hm^n(\bar{f}(N)) = 0$ or $\bar{f}(N)$ is purely $\hm^n$ unrectifiable.
\end{cor}

\end{document}